\documentclass[12pt,a4paper]{amsart}
\usepackage{amssymb,latexsym}

\usepackage{a4wide}

\theoremstyle{plain}
\newtheorem{theorem}{Theorem}[section]
\newtheorem{lemma}[theorem]{Lemma}
\newtheorem{proposition}[theorem]{Proposition}

\newtheorem*{maintheorem}{Main Theorem}

\theoremstyle{remark}

\newcommand{\C}{\ensuremath{\mathbb{C}}}
\newcommand{\R}{\ensuremath{\mathbb{R}}}
\newcommand{\E}{\ensuremath{\mathbb{E}}}
\newcommand{\g}[1]{\ensuremath{\mathfrak{#1}}}
\DeclareMathOperator{\ad}{ad}
\DeclareMathOperator{\Ad}{Ad}

\DeclareMathOperator{\Exp}{Exp}

\begin{document}
\title[Homogeneous polar foliations of complex hyperbolic
spaces]{Homogeneous polar foliations\\ of complex hyperbolic spaces}

\author[J.~Berndt]{J\"{u}rgen Berndt}
\address{Department of Mathematics\\
         King's College London\\
         United Kingdom}
\email[J\"{u}rgen Berndt]{jurgen.berndt@kcl.ac.uk}

\author[J.\,C.~D\'{\i}az-Ramos]{Jos\'{e}~Carlos D\'{\i}az-Ramos}
\address{Department of Geometry and Topology\\
         University of Santiago de Compostela\\
         Spain}
\email[Jos\'{e}~Carlos D\'{\i}az-Ramos]{josecarlos.diaz@usc.es}
\thanks{The second author has been supported by a Marie-Curie
        European Reintegration Grant (PERG04-GA-2008-239162)
        and projects MTM2009-07756
        and INCITE09207151PR (Spain)}

\subjclass[2010]{53C12, 53C35}

\keywords{Polar actions, homogeneous foliations, complex hyperbolic space}

\begin{abstract}
We prove that, up to isometric congruence, there are exactly $2n+1$
homogeneous polar foliations of the complex hyperbolic space $\C
H^n$, $n \geq 2$. We also give an explicit description of each of
these foliations.
\end{abstract}

\maketitle

\section{Introduction}

An isometric action of a connected Lie group $G$ on a Riemannian manifold $M$ is called polar if there exists a connected closed submanifold ${\mathcal S}$ of $M$ such that
\begin{itemize}
\item[(i)] ${\mathcal S}$ meets each orbit of the action, i.e. for each $p \in M$ the intersection of ${\mathcal S}$ and the orbit $G \cdot p$ of $G$ containing $p$ is nonempty;
\item[(ii)] ${\mathcal S}$ intersects each orbit orthogonally, i.e. for all $p \in {\mathcal S}$ the tangent space $T_p {\mathcal S}$ of ${\mathcal S}$ at $p$ is contained in the normal space $\nu_p (G\cdot p)$ of $G \cdot p$ at $p$.
\end{itemize}
Such a submanifold ${\mathcal S}$ is called a section of the action.
The classical paper on this topic is by Dadok \cite{Da85} who
determined all connected Lie subgroups of the orthogonal group $O(n)$
acting polarly on the Euclidean space $\E^n$. Using Dadok's result it
is not difficult to deduce a classification of all polar actions up
to orbit equivalence on the sphere $S^n$ and in the real hyperbolic
space $\R H^n$, see e.g.~\cite{Wu92}.  Podest\`{a} and Thorbergsson
classified in \cite{PT99} all polar actions on compact Riemannian
symmetric spaces of rank one up to orbit equivalence. Kollross
studied polar actions on compact Riemannian symmetric spaces of
higher rank and obtained a complete classification in \cite{Ko07} for
the case when the isometry group of the symmetric space is simple,
and in \cite{Ko09}  for the case when the symmetric space is an
exceptional compact Lie group. Recently, Kollross and Lytchak derived
the classification of polar actions on irreducible Riemannian
symmetric spaces of compact type~\cite{KoLy}. A striking difference
between the rank one case and the higher rank case is that in the
higher rank case polar actions are hyperpolar, that is, the section
${\mathcal S}$ is flat. This is not true for the rank one case.

The classification of polar actions on Riemannian symmetric spaces of
noncompact type is more complicated due to the noncompactness of the
isometry groups. One cannot expect a general correspondence using the
concept of duality between symmetric spaces of compact type and of
noncompact type. It is known that there are symmetric spaces of
compact type which do not admit any polar foliation, whereas on every
symmetric space of noncompact type there are polar foliations.
However, for actions by reductive groups Kollross established in
\cite{Ko} a correspondence. Polar actions with a fixed point on
Riemannian symmetric spaces of noncompact type have been classified
by the second author and Kollross in \cite{DK11}. Polar actions with
codimension one orbits, or equivalently, cohomogeneity one actions,
on Riemannian symmetric spaces of noncompact type are also well
understood, see \cite{BT03} and \cite{BT10}.

The motivation for this paper is to understand better the general
classification problem for the complex hyperbolic space $\C H^n =
SU(1,n)/S(U(1)U(n))$. The situation for the quaternionic hyperbolic
space $\mathbb{H}H^n$ and the Cayley hyperbolic plane $\mathbb{O}H^2$
is much more involved and the same approach does not lead to a
complete classification. We restrict here to the situation when the
orbits of the action form a foliation of $\C H^n$, in which case we
call the foliation a homogeneous polar foliation of $\C H^n$. There
are two trivial homogeneous polar foliations on $\C H^n$, one for
which the leaves are points in $\C H^n$, and one with exactly one
leaf equal to $\C H^n$. Consider a restricted root space
decomposition
$\g{g}=\g{g}_{-2\alpha}\oplus\g{g}_{-\alpha}\oplus\g{g}_0\oplus
\g{g}_\alpha\oplus\g{g}_{2\alpha}$, where $\g{g}$ is the Lie algebra
of $SU(1,n)$. The subspace $\g{g}_0$ decomposes into $\g{g}_0 =
\g{k}_0 \oplus \g{a}$, where $\g{a}$ is a one-dimensional subspace of
$\g{p}$ in a Cartan decomposition $\g{g} = \g{k} \oplus \g{p}$ and
$\g{k}_0$ is the centralizer of $\g{a}$ in $\g{k}$. The K\"{a}hler
structure of $\C H^n$ induces a complex structure on $\g{g}_\alpha$,
turning $\g{g}_\alpha$ into an $(n-1)$-dimensional complex vector
space. To each real subspace $\g{w}$ of $\g{g}_\alpha$ and each
subspace $V$ of $\g{a}$ we associate a subalgebra $\g{s}_{V,\g{w}}$
of $\g{g}$. Note that $\dim V \in \{0,1\}$ and $\dim \g{w} \in
\{0,\ldots,n-1\}$.  Let $S_{V,\g{w}}$ be the connected subgroup of
$SU(1,n)$ with Lie algebra  $\g{s}_{V,\g{w}}$. We will show that the
action of $S_{V,\g{w}}$ induces a homogeneous polar foliation of $\C
H^n$ with cohomogeneity equal to $\dim \g{w}$ if $V = \g{a}$ and
equal to $\dim \g{w} + 1$ if $V = \{0\}$. Moreover, the actions of
two such subgroups $S_{V,\g{w}}$ and $S_{V^\prime,\g{w}^\prime}$ are
orbit equivalent if and only if $\dim V = \dim V^\prime$ and $\dim
\g{w} = \dim \g{w}^\prime$. We therefore can define $S_{a,b} :=
S_{V,\g{w}}$ with $a = \dim V \in \{0,1\}$ and $b = \dim \g{w} \in
\{0,\ldots,n-1\}$, and up to isometric congruence $S_{a,b}$ is well
defined. The group $S_{0,0}$ acts transitively on $\C H^n$, and the
group $S_{1,0}$ induces a foliation of $\C H^n$ by horospheres. Our
main result states:

\begin{maintheorem}
Every nontrivial homogeneous polar foliation of $\C H^n$, $n \geq 2$, is up to isometric congruence one of the following:
\begin{itemize}
\item[(i)] The homogeneous polar foliation induced by $S_{0,b}$, $b \in \{0,\ldots,n-1\}$ - in this case the codimension of the foliation is equal to $b+1$ and all leaves are contained in horospheres of $\C H^n$;
\item[(ii)] The homogeneous polar foliation induced by $S_{1,b}$, $b \in \{1,\ldots,n-1\}$ - in this case the codimension of the foliation is equal to $b$ and no leaf is contained in a horosphere of $\C H^n$.
\end{itemize}
\end{maintheorem}

It is worthwhile to point out that none of the foliations in (i) and (ii) has a totally geodesic leaf.

This paper is organized as follows. In Section~\ref{S:CHn} we present some relevant material about the structure theory of the Lie algebra of $SU(1,n)$, and in Section~\ref{S:examples} we prove that the action of $S_{V,\g{w}}$ induces a homogeneous polar foliation of $\C H^n$. In Section~\ref{S:proof} we present the proof of the Main Theorem.

\section{Preliminaries}\label{S:CHn}

For the concepts and notation on Lie groups and Lie algebras we follow~\cite{K96}. For more information on the complex hyperbolic space and its relation to Damek-Ricci spaces see~\cite{BTV95}.

We denote by $\C H^n$, $n \geq 2$,  the complex hyperbolic space equipped with the Bergman metric normalized so that the holomorphic sectional curvature is equal to~$-1$. It can be realized as the Riemannian symmetric space $G/K$ with $G=SU(1,n)$ and $K=S(U(1)U(n))$ and a $G$-invariant Riemannian metric induced by the Killing form of the Lie algebra $\g{g}$ of $G$. Here, $G$ is the identity component of the isometry group of $\C H^n$ and $K$ is the isotropy subgroup of $G$ at a point $o\in\C H^n$, which we will fix from now on. Let $B$ denote the Killing form of $\g{g}=\g{su}(1,n)$, the Lie algebra of $SU(1,n)$. If $\g{p}$ is the orthogonal complement of $\g{k}$ in $\g{g}$ with respect to $B$, then we have a Cartan decomposition $\g{g}=\g{k}\oplus\g{p}$. The Killing form is negative definite on $\g{k}$ and positive definite on $\g{p}$. Let $\theta$ be the Cartan involution of $\g{g}$ with respect to this Cartan decomposition, that is, $\theta$ acts as the identity on $\g{k}$ and as minus the identity on $\g{p}$. Then we can define a positive definite inner product, the so-called Killing metric, on $\g{g}$ by $\langle X,Y\rangle=-B(\theta X,Y)$ for all $X$, $Y\in\g{g}$. Moreover, we have $\langle \ad(X)Y,Z\rangle=-\langle Y,\ad(\theta X)Z\rangle$ for all $X$, $Y$, $Z\in\g{g}$. As usual, $\ad$ denotes the adjoint map at Lie algebra level, $\ad(X)Y=[X,Y]$, whereas $\Ad$ will denote the adjoint map at Lie group level. Recall that $\g{p}$ can be identified with the tangent space $T_o\C H^n$ and hence can be viewed naturally as a complex vector space.

We choose a maximal abelian subspace $\g{a}$ of $\g{p}$. Then $\g{a}$ is one-dimensional. This abelian subspace induces a restricted root space decomposition $\g{g}=\g{g}_{-2\alpha}\oplus\g{g}_{-\alpha}\oplus\g{g}_0\oplus
\g{g}_\alpha\oplus\g{g}_{2\alpha}$, where $\g{g}_\lambda=\{X\in\g{g}:\ad(H)X=\lambda(H)X\text{ for all }H\in\g{a}\}$ and $\lambda \in \g{a}^*$. Recall that $[\g{g}_\lambda,\g{g}_\mu]=\g{g}_{\lambda+\mu}$. The set $\Sigma=\{-2\alpha,-\alpha,\alpha,2\alpha\}$ is called the set of roots. It is known that $\theta\g{g}_\lambda=\g{g}_{-\lambda}$ for any $\lambda\in\Sigma\cup\{0\}$, and $\g{g}_0=\g{k}_0\oplus\g{a}$, where $\g{k}_0=\g{g}_0\cap\g{k}$. It can be seen that $\g{k}_0$ is isomorphic to $\g{u}(n-1)$, and that $\g{g}_{2\alpha}$ is one-dimensional. From now on, we introduce an ordering in $\Sigma$ such that $\alpha$ and $2\alpha$ are positive roots. This choice selects precisely one unit vector $B$ in $\g{a}$ for which $\alpha(B)=1/2$.

Let $\g{n}=\g{g}_\alpha\oplus\g{g}_{2\alpha}$, which is a nilpotent subalgebra of $\g{g}$ isomorphic to the $(2n-1)$-dimensional Heisenberg algebra. Then $\g{g}=\g{k}\oplus\g{a}\oplus\g{n}$ is an Iwasawa decomposition of $\g{g}$. It is known that the connected subgroup $AN$ of $G$ whose Lie algebra is $\g{a}\oplus\g{n}$ acts simply transitively on $\C H^n$. We endow $AN$, and hence $\g{a}\oplus\g{n}$, with the left-invariant metric $\langle\,\cdot\,,\,\cdot\,\rangle_{AN}$ and the complex structure $J$ that make $\C H^n$ and $AN$ isometric. If $X$, $Y\in\g{a}\oplus\g{n}\cong T_{1}AN$ are considered as left-invariant vector fields then the relation between the Killing form and the inner product on $\g{a}\oplus\g{n}$ is given by $\langle{X},{Y}\rangle_{AN}
=\langle X_{\g{a}},Y_{\g{a}}\rangle
+\frac{1}{2}\langle X_{\g{n}},Y_{\g{n}}\rangle$,
where subscript means orthogonal projection. The complex structure $J$ on $\g{a}\oplus\g{n}$ is such that $\g{g}_\alpha$ is $J$-invariant, thus $\g{g}_\alpha$ is a complex vector subspace, and $J\g{a}=\g{g}_{2\alpha}$. Therefore we can define $Z=JB\in\g{g}_{2\alpha}$. Notice that, while $B$ and $Z$ are unit vectors with respect to the $\g{a}\oplus\g{n}$ metric, we have $\langle B,B\rangle=1$ but $\langle Z,Z\rangle=2$ with respect to the Killing metric. With this notation, the Lie bracket in $\g{a}\oplus\g{n}$ can be written as
\[
[aB+U+xZ,bB+V+yZ]=
-\frac{b}{2}U+\frac{a}{2}V
+\left(-bx+ay+\frac{1}{2}\langle JU,V\rangle\right)Z,
\]
where $a$, $b$, $x$, $y\in\R$, and $U$, $V\in\g{g}_\alpha$.

We can also define $\g{p}_\lambda=(1-\theta)\g{g}_\lambda\subset\g{p}$ and $\g{k}_\lambda=(1+\theta)\g{g}_\lambda\subset\g{k}$. Then we have  $\g{p}=\g{a}\oplus\g{p}_\alpha\oplus\g{p}_{2\alpha}$, $\g{p}_\alpha$ is complex, and $\g{p}_{2\alpha}$ is one-dimensional. If $i$ denotes the complex structure of $\g{p}$ then we have
\begin{equation}\label{E:J}
iB=\frac{1}{2}(1-\theta)Z\quad\text{ and }\quad i(1-\theta)U=(1-\theta)JU.
\end{equation}

There are exactly two conjugacy classes of Cartan subalgebras in $\g{g}$. We are interested in the so-called maximally noncompact ones. These are, up to conjugacy, of the form $\g{t}\oplus\g{a}$, where $\g{t}$ is some abelian subspace of $\g{k}_0$. Then, $[\g{t},\g{a}]=[\g{t},\g{g}_{2\alpha}]=0$, and $[\g{t},\g{g}_\alpha]\subset\g{g}_\alpha$. The vector space
$\g{t}\oplus\g{a}\oplus\g{n}$ turns out to be a maximal solvable subalgebra of $\g{g}$, of the so-called maximal noncompact type.

To finish this section, we prove the following result, which will be used several times henceforth.

\begin{lemma}
We have
\begin{align}
&[\theta U,Z]=-JU\text{ for all }U\in\g{g}_\alpha;
\label{E:bracket theta X, Z}\\[1ex]
&\langle T,(1+\theta)[\theta U,V]\rangle=2\langle[T,U],V\rangle\text{ for all $T\in\g{t}$ and $U$, $V\in\g{g}_\alpha$}.\label{E:auxiliary}
\end{align}
\end{lemma}

\begin{proof}
Assume that $U\in\g{g}_\alpha$. It is clear that $[\theta
U,Z]\in\g{g}_\alpha$. Moreover, if $V$ is another element in
$\g{g}_\alpha$ we have
$\langle[\theta U,Z],V\rangle=-\langle Z,[U,V]\rangle
=-\frac{1}{2}\langle JU,V\rangle\langle Z,Z\rangle=-\langle JU,V\rangle$,
where we have used $\langle Z,Z\rangle=2$. This
implies~(\ref{E:bracket theta X, Z}).

Now let $T\in\g{t}$ and $U$, $V\in\g{g}_\alpha$. Then
\[
\langle T,(1+\theta)[\theta U,V]\rangle
=\langle T,[\theta U,V]\rangle+\langle T,[U,\theta V]\rangle
=-\langle[U,T],V\rangle+\langle[V,T],U\rangle
=2\langle[T,U],V\rangle,
\]
from where~(\ref{E:auxiliary}) follows.
\end{proof}

\section{Polar foliations}\label{S:examples}

The following criterion for homogeneous polar foliations was proved in \cite{BDT10}.

\begin{theorem}\label{T:polarFoliations}
Let $M = G/K$ be a Riemannian symmetric space of noncompact type and $H$ be a
connected closed subgroup of $G$ whose orbits form a homogeneous foliation
$\mathcal{F}$ of $M$. Consider the corresponding Cartan decomposition $\g{g} =
\g{k} \oplus \g{p}$ and define
\[
{\g{h}_{\g{p}}^\perp}=\{\,\xi\in\g{p}:\langle\xi,Y\rangle
=0\mbox{ \normalfont for
all }Y\in\g{h}\,\}.
\]
Then the action of $H$ on $M$ is polar if and only if
${\g{h}_{\g{p}}^\perp}$ is a Lie triple system in $\g{p}$ and
$\g{h}$ is orthogonal to the subalgebra
$[{\g{h}_{\g{p}}^\perp},{\g{h}_{\g{p}}^\perp}]\oplus{\g{h}_{\g{p}}^\perp}$
of $\g{g}$. In this case, let ${H_{\g{p}}^\perp}$ be the
connected subgroup of  $G$ with Lie algebra
$[{\g{h}_{\g{p}}^\perp},{\g{h}_{\g{p}}^\perp}]
\oplus{\g{h}_{\g{p}}^\perp}$. Then the orbit ${\mathcal S} =
{H_{\g{p}}^\perp} \cdot o$ is a section of the $H$-action on
$M$.
\end{theorem}

As a first application of this theorem we construct examples of
homogeneous polar foliations of the complex hyperbolic space $\C
H^n$. Here and henceforth the symbol $\ominus$ denotes orthogonal complement.

\begin{theorem}\label{T:examples}
Let $\g{w}$ be a real subspace of $\g{g}_\alpha$ and $V$ be a subspace
of $\g{a}$ (then, either $V=\{0\}$ or $V=\g{a}$). Let $S_{V,\g{w}}$ be
the connected subgroup of $AN$ whose Lie algebra is
$\g{s}_{V,\g{w}}=V\oplus(\g{g}_\alpha\ominus\g{w})\oplus\g{g}_{2\alpha}$.
Then $S_{V,\g{w}}$ induces a homogeneous polar foliation of $\C H^n$.
\end{theorem}

\begin{proof}
It is clear that $\g{s}_{V,\g{w}}$ is a subalgebra of $\g{a} \oplus \g{n}$. Since $S_{V,\g{w}}$ is contained in $AN$ it is obvious that it
induces a homogeneous foliation of $\C H^n$. So, it remains to use
Theorem~\ref{T:polarFoliations} to show that it also acts polarly on
$\C H^n$.

It is easy to calculate that
\[
(\g{s}_{V,\g{w}})_{\g{p}}^\perp=(\g{a}\ominus V)\oplus(1-\theta)\g{w}.
\]
The subspace above is real in $\g{p}$ (use for example~\eqref{E:J}), and hence it is a Lie triple system of
$\g{p}$. Now, if $W$, $\tilde{W}\in\g{w}$, we have
$[(1-\theta)W,(1-\theta)\tilde{W}]=(1+\theta)[W,\tilde{W}]-(1+\theta)[\theta
W,\tilde{W}]=-(1+\theta)[\theta W,\tilde{W}]\in\g{k}_0$, which is
orthogonal to $\g{a}\oplus\g{n}$ and thus to $\g{s}_{V,\g{w}}$. Also,
$[B,(1-\theta)W]=\frac{1}{2}(1+\theta)W$, and since $\g{w}$ is
orthogonal to $\g{s}_{V,\g{w}}$, so is the latter element of
$\g{k}_\alpha$. This finishes the proof of the theorem.
\end{proof}

Denote by $K_0$ the subgroup of $K$ with Lie algebra $\g{k}_0$. The adjoint action of $K_0$ on $\g{g}_\alpha$ is isomorphic to the standard action of $U(n-1)$ on $\C^{n-1}$. This implies that for any two real subspaces $\g{w}$ and $\g{w}^\prime$ of $\g{g}_\alpha$ with $\dim \g{w} = \dim \g{w}^\prime$ the subgroups $S_{V,\g{w}}$ and $S_{V,\g{w}^\prime}$ are conjugate in $SU(1,n)$. Hence the corresponding two homogeneous polar foliations on $\C H^n$ are isometrically congruent to each other.

We will now proceed with the proof of the Main Theorem.

\section{Proof of the Main Theorem}\label{S:proof}

Assume that $H$ is a connected closed subgroup of $G=SU(1,n)$ acting polarly on $\C H^n$
in such a way that the orbits of $H$ form a homogeneous foliation of $\C H^n$.
We denote by $\g{h}$ the Lie algebra of $H$. According to \cite{BDT10} we have:

\begin{theorem}\label{T:ReduceToSolvable}
Let $M$ be a Hadamard manifold and let $H$ be a connected closed
subgroup of $I(M)$ whose orbits form a homogeneous foliation
${\mathcal F}$ of $M$. Then there exists a connected closed solvable
subgroup $S$ of $H$ such that the leaves of ${\mathcal F}$ coincide
with the orbits of $S$. Moreover, all the orbits of $S$ or of $H$ are
principal.
\end{theorem}

Theorem~\ref{T:ReduceToSolvable} implies that we may assume that $\g{h}$ is solvable. Thus $\g{h}$ is contained in a maximal solvable subalgebra of $\g{g}$. We can give a more precise statement in our situation.

\begin{proposition}\label{P:Borel}
Let $M$ be a symmetric space on noncompact type and rank one. Assume that $H$ acts polarly on $M$ inducing a foliation. Then the action of $H$ is orbit equivalent to the action of a Lie group whose Lie algebra is contained in a maximally noncompact Borel subalgebra.
\end{proposition}

\begin{proof}
We write $M=G/K$, where $G$ is the connected component of the isometry group of $M$. By Theorem~\ref{T:ReduceToSolvable} we may assume that $H$ is solvable and closed in $G$.

Let $\g{b}$ be a maximal solvable subalgebra of $\g{g}$
containing $\g{h}$. Then there exists a Cartan decomposition
$\g{g}=\g{k}\oplus\g{p}$ such that
$\g{b}=\g{t}\oplus\tilde{\g{a}}\oplus\tilde{\g{n}}$, with
$\g{t}\subset\g{k}$, $\tilde{\g{a}}\subset\g{p}$, and such that
$\g{t}\oplus\tilde{\g{a}}$ is Cartan subalgebra of $\g{g}$
\cite{M61}. Here, $\tilde{\g{n}}$ is defined as $\tilde{\g{n}}
=\bigoplus_{\tilde{\lambda}\in\tilde{\Sigma}^+}
\tilde{\g{g}}_{\tilde{\lambda}}$, where $\tilde{\Sigma}$ is the set of  roots with respect to $\tilde{\g{a}}$, $\tilde{\g{g}}_{\tilde{\lambda}}=\{X\in\g{g}:
\ad(H)X=\tilde{\lambda}(H)X\text{ for all }H\in\tilde{\g{a}}\}$, and
$\tilde{\Sigma}^+$ is the set of positive roots with respect to some choice  of ordering. Since $M$ has rank one, either $\tilde{\g{a}}=\{0\}$ or
$\dim\tilde{\g{a}}=1$. If $\tilde{\g{a}}=\{0\}$ then $\g{b}=\g{t}$ is
compact. Hence, the action of $H$ is orbit equivalent to the action
of a compact subgroup, and thus has a fixed point by Cartan's fixed
point theorem. Since $H$ induces a foliation, it follows that the orbits of $H$ are points. Therefore, we may assume
$\dim\tilde{\g{a}}=1$, and thus $\g{b}$ is maximally noncompact. We
conclude that $\g{h}$ is a subalgebra of
$\g{t}\oplus\g{a}\oplus\g{n}$, where $\g{a}\subset\g{p}$ is a maximal
abelian subspace.
\end{proof}

In view of Proposition~\ref{P:Borel} we may assume that $\g{h}$ is a subalgebra of
$\g{t}\oplus\g{a}\oplus\g{n}$, where $\g{a}\subset\g{p}$ is a maximal
abelian subspace and $\g{t}\oplus\g{a}$ is a maximally noncompact Cartan subalgebra of $\g{g}$.

Let $\mathcal{S}$ be a section of the action of $H$. Then, according to Theorem~\ref{T:polarFoliations}, $\mathcal{S}=H_{\g{p}}^\perp\cdot o$, where $H_{\g{p}}^\perp$ is the connected subgroup of $G=SU(1,n)$ whose Lie algebra is $[\g{h}_{\g{p}}^\perp,\g{h}_{\g{p}}^\perp] \oplus \g{h}_{\g{p}}^\perp \subset \g{k} \oplus \g{p} = \g{g}$, where $\g{h}_{\g{p}}^\perp=\{\,\xi\in\g{p}:\langle\xi,Y\rangle
=0\mbox{ \normalfont for
all }Y\in\g{h}\,\}\subset\g{p}$. Since $\g{h}_{\g{p}}^\perp$ is a Lie triple system it follows that
$\g{h}_{\g{p}}^\perp$ is either real or complex, and thus $\mathcal{S}$ is also either real or complex. We deal first with the complex case.

\begin{proposition}
If the section $\mathcal{S}$ of $H$ is complex, then either $H$ acts transitively on $\C H^n$ or the orbits of $H$ are points.
\end{proposition}

\begin{proof}
Let $\pi\colon\g{h}_{\g{p}}^\perp\to\g{a}\oplus\g{p}_{2\alpha}$ be the orthogonal projection of $\g{h}_{\g{p}}^\perp$ onto $\g{a}\oplus\g{p}_{2\alpha}$. Note that $\pi$ is $J$-linear, and thus $\pi(\g{h}_{\g{p}}^\perp)$ is either $\{0\}$ or $\g{a}\oplus\g{p}_{2\alpha}$.

Assume first that $\pi(\g{h}_{\g{p}}^\perp)=\{0\}$. Then it is clear that $\g{h}_{\g{p}}^\perp$ is a complex subspace of $\g{p}_\alpha$. Since $B$ and $Z$ are orthogonal to $\g{h}_{\g{p}}^\perp$, there exist $S$, $T\in\g{t}$ such that $S+B$,
$T+Z\in\g{h}$. Then,
$Z=[S+B,T+Z]\in\g{h}$. Let $V\in\g{g}_{\alpha}$ such that
$(1-\theta)V\in\g{h}_{\g{p}}^\perp$. Theorem~\ref{T:polarFoliations}
implies
\[
0=\langle Z,[(1-\theta)V,(1-\theta)JV]\rangle
=\langle Z,[V,JV]\rangle=\frac{1}{2}\langle V,V\rangle\langle Z,Z\rangle
=\lVert V\rVert^2,
\]
and so $\g{h}_{\g{p}}^\perp=\{0\}$. Hence the section is zero-dimensional and therefore the action of $H$
must be transitive.

Now assume that $\pi(\g{h}_{\g{p}}^\perp)=\g{a}\oplus\g{p}_{2\alpha}$. As $\pi$ is $J$-linear, $\ker\pi$ is a complex subspace of $\g{p}_\alpha$. Then $\g{h}_{\g{p}}^\perp\ominus\ker\pi$, the orthogonal complement of $\ker\pi$ in $\g{h}_{\g{p}}^\perp$, has complex dimension one, so there is a unique vector $Y\in\g{h}_{\g{p}}^\perp\ominus\ker\pi$ such that $\pi(Y)=B$. Thus, we can write $Y=B+(1-\theta)X$, with $X\in\g{g}_\alpha$ and $X$ (or $(1-\theta)X$) orthogonal to $\ker\pi$. In view of~\eqref{E:J}, we can also write
\[
\g{h}_{\g{p}}^\perp=\C(B+X)\oplus\ker\pi
=\R(B+(1-\theta)X)\oplus\R((1-\theta)JX
+\frac{1}{2}(1-\theta)Z)\oplus\ker\pi.
\]

It is obvious that
$\lVert X\rVert^2 B-X$, and $-JX+\lVert X \rVert^2 Z$ are
orthogonal  to $\g{h}_{\g{p}}^\perp$ (recall that
$\lVert B\rVert^2=1$, and $\lVert Z\rVert^2=2$, where
$\lVert\cdot\rVert$ is the norm induced by
$\langle\,\cdot\,,\,\cdot\,\rangle$). Hence, there exist $S$,
$T\in\g{t}$ such that $S+\lVert X\rVert^2 B-X$, $T-
JX+\lVert X \rVert^2 Z\in\g{h}$.

Using~\eqref{E:bracket theta X, Z} we get
\begin{align*}
&[\g{h}_{\g{p}}^\perp,\g{h}_{\g{p}}^\perp]\ni
\left[B+(1-\theta)X,(1-\theta)JX+\frac{1}{2}(1-\theta)Z\right]\\
&\quad=\frac{1}{2}(1+\theta)JX+\frac{1}{2}(1+\theta)Z+(1+\theta)[X,JX]
-(1+\theta)[\theta X,JX]-\frac{1}{2}(1+\theta)[\theta X,Z]\\
&\quad=-(1+\theta)[\theta X,JX]+(1+\theta)JX
+\frac{1}{2}\left(1+{\lVert X\rVert^2}\right)(1+\theta)Z.
\end{align*}

Taking inner product of the previous vector with $T-JX+\lVert X
\rVert^2 Z\in\g{h}$, and using the fact that $\ad(T)$ is
skew-symmetric, and~\eqref{E:auxiliary}, we get by
Theorem~\ref{T:polarFoliations}
\[
0=-\langle (1+\theta)[\theta X,JX],T\rangle
-{\lVert X\rVert^2}
+\frac{\lVert X\rVert^2}{2}
\left(1+\lVert X\rVert^2\right)\lVert Z\rVert^2
=-2\langle [T,X],JX\rangle
+{\lVert X\rVert^4}
\]
Therefore we obtain $\langle [T,X],JX\rangle=\lVert
X\rVert^4/2$.

On the other hand, since $\g{h}$ is a Lie algebra, and
$[\g{t},\g{g}_{2\alpha}]=0$, we have that
\begin{align*}
\g{h}\ni{} &[S+\lVert X\rVert^2 B-X,T-JX+\lVert X \rVert^2 Z]\\
&\quad=[T,X]-[S,JX]-\frac{\lVert X\rVert^2}{2}JX
+\lVert X\rVert^2\left(\lVert X\rVert^2+\frac{1}{2}\right)Z.
\end{align*}
Taking inner product with
$(1-\theta)JX+\frac{1}{2}(1-\theta)Z\in\g{h}_{\g{p}}^\perp$,
and using the fact that $\ad(S)$ is skew-symmetric, yields
\[
0=\langle [T,X],JX\rangle-\frac{\lVert X\rVert^4}{2}
+\lVert X\rVert^2\left(\lVert X\rVert^2+\frac{1}{2}\right)
=\lVert X\rVert^2\left(\lVert X\rVert^2+\frac{1}{2}\right).
\]
Therefore $X=0$, and $\g{h}_{\g{p}}^\perp=\g{a}\oplus\g{p}_{2\alpha}\oplus\ker\pi$. Since $\ker\pi$ is complex we can find
$S$, $T\in\g{t}$, and $V\in\g{g}_{\alpha}$, such that $S+V$,
$T+JV\in\g{h}$. Then
\begin{align*}
0=\langle[S+V,T+JV],(1-\theta)Z\rangle
=\langle [S,JV]-[T,V]+\frac{1}{2}\lVert V\rVert^2 Z,
(1-\theta)Z\rangle
=\lVert V\rVert^2.
\end{align*}
This implies $\g{w}=\g{p}_{\alpha}$ and $\g{h}\subset\g{t}$.
Therefore $H$ is compact. By Cartan's fixed point theorem $H$ has a
fixed point and thus fixes $\C H^n$ pointwise.
\end{proof}

Altogether this proves that a homogeneous polar foliation of
the complex hyperbolic space cannot have a nontrivial complex
section. \medskip

From now on we assume that the section is real, that is, $\g{h}_{\g{p}}^\perp$ is a real subspace of $\g{p}$.
Consider again the orthogonal projection $\pi\colon\g{h}_{\g{p}}^\perp\to\g{a}\oplus\g{p}_{2\alpha}$. Two different possibilities arise.

\subsubsection*{Case 1: $\dim\pi(\g{h}_{\g{p}}^\perp)=2$}\hfil
We will show that this possibility cannot hold. We start with an algebraic result.

\begin{lemma}\label{L:Real1}
The subspace $\g{h}_{\g{p}}^\perp$ can be written as
\[\g{h}_{\g{p}}^\perp=\R(B+(1-\theta)X)\oplus\R((1-\theta)Y+(1-\theta)Z)
\oplus(1-\theta)\g{w},
\] where $X$, $Y\in\g{g}_\alpha$, $\g{w}$ is a real subspace of $\g{g}_\alpha$, $\C X$ and $\C Y$ are perpendicular to $\g{w}$, and $\langle X,JY\rangle=1$.
\end{lemma}

\begin{proof}
Since $\ker\pi\subset\g{p}_\alpha$, we can find a subspace $\g{w}\subset\g{g}_\alpha$ such that $\ker\pi=(1-\theta)\g{w}$.
By hypothesis, there exists a unique $\xi\in\g{h}_{\g{p}}^\perp\ominus\ker\pi$ such that $\pi(\xi)=B$. Hence, there exists $X\in\g{g}_\alpha$ such that $\xi=B+(1-\theta)X\in\g{h}_{\g{p}}^\perp$. Similarly, there exists a unique $\eta\in\g{h}_{\g{p}}^\perp\ominus\ker\pi$ such that $\pi(\eta)=(1-\theta)Z$, and thus, there exists $Y\in\g{g}_\alpha$ such that $\eta=(1-\theta)Z+(1-\theta)Y\in\g{h}_{\g{p}}^\perp$. Clearly, $\g{h}_{\g{p}}^\perp=\R\xi\oplus\R\eta\oplus\ker\pi
=\R(B+(1-\theta)X)\oplus\R((1-\theta)Y+(1-\theta)Z)
\oplus(1-\theta)\g{w}$.

Since $\g{h}_{\g{p}}^\perp$ is real, using~\eqref{E:J} we get $0=\langle\xi,J\eta\rangle=-2+2\langle X,JY\rangle$. Let $W\in\g{w}\subset\g{g}_\alpha$ be arbitrary. By definition of $\xi$, we get $0=\langle (1-\theta)W,\xi\rangle=2\langle W,X\rangle$, and since $\g{h}_{\g{p}}^\perp$ is real, $0=\langle (1-\theta)W,J\xi\rangle=2\langle W,JX\rangle$. This implies that $\C X$ is orthogonal to $\g{w}$. We can prove in a similar way that $\C Y$ is orthogonal to $\g{w}$, and hence the result follows.
\end{proof}

In view of Lemma~\ref{L:Real1} our first observation is that
\begin{equation}\label{E:orthogonal to hpperp}
\text{for each }U\in\g{g}_\alpha\ominus\g{w},\quad -\langle U,X\rangle B+U-\frac{1}{2}\langle U,Y\rangle Z\text{ is orthogonal to $\g{h}_{\g{p}}^\perp$}.
\end{equation}
Hence, for each
$U\in\g{g}_\alpha\ominus\g{w}$ there exists $T_U\in\g{t}$ such that
$T_U-\langle U,X\rangle B+U-\frac{1}{2}\langle U,Y\rangle
Z\in\g{h}$.

Using~\eqref{E:bracket theta X, Z} we also get
\begin{equation}\label{E:bracket}
[B+(1-\theta)X,(1-\theta)Y+(1-\theta)Z]=
-(1+\theta)[\theta X,Y]+\frac{1}{2}(1+\theta)Y+(1+\theta)JX+\frac{1}{2}(1+\theta)Z.
\end{equation}

\begin{lemma}\label{L:Y}
We have that $Y\in\C X$. More explicitly, there exists $\gamma\in\R$ such that
\[
Y=\gamma X-\frac{1}{\lVert X\rVert^2}JX.
\]
\end{lemma}

\begin{proof}
According to~\eqref{E:orthogonal to hpperp}, there exist $T_X$, $T_Y\in\g{t}$ such that
\[
T_X-{\lVert X\rVert^2}B+X-\frac{1}{2}\langle X,Y\rangle Z,\quad
T_Y-\langle X,Y\rangle B+Y-\frac{\lVert Y\rVert^2}{2}Z\in\g{h}.
\]

Since $\g{h}$ is orthogonal to
$[\,\g{h}_{\g{p}}^\perp,\g{h}_{\g{p}}^\perp]$ by
Theorem~\ref{T:polarFoliations}, taking inner product of
$T_Y-\langle X,Y\rangle B+Y-\frac{1}{2}\lVert Y\rVert Z$
with~\eqref{E:bracket}, and using~\eqref{E:auxiliary} we
get $0=-\langle T_Y,(1+\theta)[\theta X,Y]\rangle
=-2\langle[T_Y,X],Y\rangle{-1}$, so $\langle[T_Y,Y],X\rangle=1/2$.

Now, since $\g{h}$ is a Lie algebra
\begin{align*}
\g{h}\ni{}&[T_X-\lVert X\rVert^2 B+X-\frac{1}{2}\langle X,Y\rangle Z,
T_Y-\langle X,Y\rangle B+Y-\frac{\lVert Y\rVert^2}{2}Z]\\
&{}=[T_X,Y]-[T_Y,X]+\frac{1}{2}(\langle X,Y\rangle X-\lVert X\rVert^2 Y)
+\left(\frac{\lVert X\rVert^2\lVert Y\rVert^2}{2}
-\frac{\langle X,Y\rangle^2}{2}-\frac{1}{2}\right)Z.
\end{align*}
Taking inner product with
$(1-\theta)Y+(1-\theta)Z\in\g{h}_{\g{p}}^\perp$, using the
skew-symmetry of $\ad(T_X)$ and $\langle[T_Y,Y],X\rangle=1/2$, we get
\[
0=\langle [T_X,Y]-[T_Y,X],Y\rangle+\frac{1}{2}(\lVert X\rVert^2\lVert Y\rVert^2
-\langle X,Y\rangle^2)+\langle JX,Y\rangle
=\frac{1}{2}(\lVert X\rVert^2\lVert Y\rVert^2
-\langle X,Y\rangle^2-1).
\]
This implies
\[
\lVert X\rVert^2\lVert Y\rVert^2-\langle X,Y\rangle^2=1=\langle JX,Y\rangle^2.
\]
We may write $Y=\gamma X+\delta JX+E$ with $\gamma$, $\delta\in\R$
and $E\in\g{g}_\alpha\ominus\C X$. Then, the equation above reads
$\lVert X\rVert^2\lVert E\rVert^2=0$, and since $X\neq 0$ we get
$E=0$. As $\langle X,JY\rangle=1$ this readily implies
the result.
\end{proof}

Let $g\in G$ be an arbitrary isometry. The groups $H$ and
$I_g(H)=gHg^{-1}$ are conjugate, and since $H$ induces a foliation so
does $I_g(H)$. Since all the orbits of a homogeneous foliation
are principal and $H$ and $I_g(H)$ are conjugate, so are their
isotropy groups at the origin, $H_o$ and $I_g(H)_o$ respectively. In
particular, their Lie algebras have the same dimension, that
is, $\dim(\g{h}\cap\g{k})=\dim(\Ad(g)(\g{h})\cap\g{k})$. From
now on let us choose
\[
g=\Exp\left(-\frac{2}{\lVert X\rVert^2}X\right).
\]
We will see that
$\dim(\g{h}\cap\g{k})<\dim(\Ad(g)(\g{h})\cap\g{k})$, thus
leading to a contradiction.

\begin{lemma}\label{L:ht}
We have $[\g{h}\cap\g{t},X]=0$. In particular, $\Ad(g)(\g{h}\cap\g{t})=\g{h}\cap\g{t}$.
\end{lemma}

\begin{proof}
Let $S\in\g{h}\cap\g{t}$. We will indeed show $[S,X]=0$, which clearly implies $\Ad(g)(S)=S$.

Take $U\in\g{g}_\alpha\ominus(\C
X\oplus\g{w})$, and according to~\eqref{E:orthogonal to hpperp}, let $T_U\in\g{t}$ be such that
$T_U+U\in\g{h}$. Then, $[S,U]=[S,T_U+U]\in\g{h}$ and thus
$0=\langle[S,U],B+(1-\theta)X\rangle=-\langle[S,X],U\rangle$
by the skew-symmetry of $\ad(S)$. This implies that $[S,X]\in\R
JX\oplus\g{w}$.

Since $\g{h}$ is orthogonal to
$[\g{h}_{\g{p}}^\perp,\g{h}_{\g{p}}^\perp]$ we have, applying~\eqref{E:auxiliary}, \eqref{E:bracket}, and Lemma~\ref{L:Y}
\begin{align*}
0   &=\langle S,[B+(1-\theta)X,(1-\theta)Y+(1-\theta)Z]\rangle\\
    &=-\langle S,(1+\theta)[\theta X,Y]\rangle
    =-2\langle [S,X],Y\rangle
    =\frac{2}{\lVert X\rVert^2}\langle [S,X],JX\rangle,\text{ and}\\
0   &=\langle S,[B+(1-\theta)X,(1-\theta)W]\rangle
    =\langle S,-(1+\theta)[\theta X,W]+\frac{1}{2}(1+\theta)W\rangle
    =-2\langle [S,X],W\rangle,
\end{align*}
for any $W\in\g{w}$. Altogether this implies $[S,X]=0$, from where we get the result.
\end{proof}

It follows from Lemma~\ref{L:ht} that $\g{h}\cap\g{t}=\Ad(g)(\g{h}\cap\g{t})\subset\Ad(g)(\g{h})\cap\g{k}$, so $\dim(\g{h}\cap\g{k})=\dim(\g{h}\cap\g{t})\leq\dim(\Ad(g)(\g{h})\cap\g{t})
\leq\dim(\Ad(g)(\g{h})\cap\g{k})$. From~\eqref{E:orthogonal to hpperp}, there exists $T_{JX}\in\g{t}$ such that
$T_{JX}+JX+\frac{1}{2}Z\in\g{h}$. We will show the following two facts:
\[
\Ad(g)(T_{JX}+JX+\frac{1}{2}Z)=T_{JX}\in\Ad(g)(\g{h})\cap\g{t},\quad\text{ and }\quad
T_{JX}\notin\g{h}\cap\g{t}.
\]
This, and the previous inequalities, exhibit the fact that $\dim(\g{h}\cap\g{k})<\dim(\Ad(g)(\g{h})\cap\g{k})$ which gives the desired
contradiction.

We have by~\eqref{E:auxiliary} and~\eqref{E:bracket}
\begin{align*}
0   &=\langle T_{JX}+JX+\frac{1}{2}Z,[B+(1-\theta)X,(1-\theta)Y
        +(1-\theta)Z]\rangle\\
    &=-2\langle [T_{JX},X],Y\rangle+\lVert X\rVert^2
    =\frac{2}{\lVert X\rVert^2}\langle [T_{JX},X],JX\rangle+\lVert X\rVert^2,
\end{align*}
and hence,
\[
\langle [T_{JX},X],JX\rangle=-\frac{\lVert X\rVert^4}{2}.
\]
In particular, since $[\g{h}\cap\g{t},X]=0$ by Lemma~\ref{L:ht} and $X\neq 0$, this implies $T_{JX}\notin\g{h}\cap\g{t}$.

Let us then show $\Ad(g)(T_{JX}+JX+\frac{1}{2}Z)=T_{JX}$. Let
$U\in\g{g}_\alpha\ominus(\C X\oplus \g{w})$ and by~\eqref{E:orthogonal to hpperp} take $T_U\in\g{t}$
such that $T_U+U\in\g{h}$. Then, by~\eqref{E:auxiliary} and~\eqref{E:bracket},
\[
0   =\langle T_{U}+U,[B+(1-\theta)X,(1-\theta)Y+(1-\theta)Z]\rangle\\
    =-2\langle[T_U,X],Y\rangle
    =\frac{2}{\lVert X\rVert^2}\langle [T_U,X],JX\rangle,
\]
and thus $\langle [T_U,X],JX\rangle=0$. On the other hand,
\[[T_{JX},U]-[T_U,JX]=[T_{JX}+JX+\frac{1}{2}Z,T_U+U]\in\g{h},\] and
hence
\[0=\langle[T_{JX},U]-[T_U,JX],B+(1-\theta)X\rangle
=\langle[T_{JX},U],X\rangle
-\langle[T_U,JX],X\rangle=-\langle[T_{JX},X],U\rangle,\] thus showing
$[T_{JX},X]\in \R JX\oplus\g{w}$. Now, for any $W\in\g{w}$,
\begin{align*}
0   &=\langle T_{JX}+JX+\frac{1}{2}Z,[B+(1-\theta)X,(1-\theta)W]\rangle\\
    &=-\langle T_{JX}+JX+\frac{1}{2}Z,
        -(1+\theta)[\theta X,W]+\frac{1}{2}(1+\theta)W\rangle\\
    &=-\langle T_{JX},(1+\theta)[\theta X,W]\rangle
    =-2\langle [T_{JX},X],W\rangle,
\end{align*}
using~\eqref{E:auxiliary} once again. This, together with $\langle
[T_{JX},X],JX\rangle=-\frac{\lVert X\rVert^4}{2}$ implies
\[
[T_{JX},X]=-\frac{\lVert X\rVert^2}{2}JX.
\]
Now, taking this into account, and $[X,JX]=\frac{\lVert
X\rVert^2}{2}Z$, we obtain
\begin{align*}
&\Ad(g)(T_{JX}+JX+\frac{1}{2}Z)
=e^{\ad\left(\frac{-2}{\lVert X\rVert^2}X\right)}(T_{JX}+JX+\frac{1}{2}Z)\\
&\quad=T_{JX}+JX+\frac{1}{2}Z-\frac{2}{\lVert X\rVert^2}([X,T_{JX}]+[X,JX])
+\frac{2}{\lVert X\rVert^4}[X,[X,T_{JX}]]=T_{JX},
\end{align*}
which leads to the desired contradiction.

Altogether we have proved

\begin{proposition}
$\dim\pi(\g{h}_{\g{p}}^\perp)\leq 1$.
\end{proposition}

\subsubsection*{Case 2: $\dim\pi(\g{h}_{\g{p}}^\perp)\leq 1$}
In this case, there exists a nonzero vector $\xi\in\g{h}_{\g{p}}^\perp\ominus\ker\pi$, which we can write as $\xi=aB+(1-\theta)X+b(1-\theta)Z$, with $X\in\g{g}_\alpha$, $X$ orthogonal to $\ker\pi$, and $a$, $b\in\R$. Clearly, $a=b=0$ if $\pi(\g{h}_{\g{p}}^\perp)=0$; otherwise, at least one of these numbers is nonzero. Let $\g{w}$ be a subspace of $\g{g}_\alpha$ such that $\ker\pi=(1-\theta)\g{w}$. Note that, since $\g{h}_{\g{p}}^\perp$ is real, $\g{w}$ is also real, and $0=\langle J\xi,(1-\theta)W\rangle=2\langle JX,W\rangle$ for any $W\in\g{w}$, so $\R X\oplus\g{w}$ is a real subspace of $\g{g}_\alpha$.

\begin{lemma}\label{L:b=0}
We have that $\g{h}_{\g{p}}^\perp=\R(aB+(1-\theta)X)\oplus(1-\theta)\g{w}$, where $a\in\R$, $X\in\g{g}_\alpha\ominus\g{w}$, and $\R X\oplus\g{w}$ is a real subspace of $\g{g}_\alpha$.
\end{lemma}

\begin{proof}
We only have to prove that $b=0$. Assume, on the contrary, that $b\neq 0$.
It is easy to see that $\lVert X\rVert^2 B-aX$ and $-2bX+\lVert
X\rVert^2 Z$ are orthogonal to $\g{h}_{\g{p}}^\perp$, so there
exist $S$, $T\in\g{t}$ such that $S+\lVert X\rVert^2 B-aX$,
$T-2bX+\lVert X\rVert^2 Z\in\g{h}$. Then
\[
\g{h}\ni[S+\lVert X\rVert^2 B-aX,T-2bX+\lVert X\rVert^2 Z]
=-2b[S,X]+a[T,X]-b\lVert X\rVert^2 X+\lVert X\rVert^4 Z.
\]
Taking inner product of this vector with
$aB+(1-\theta)X+b(1-\theta)Z\in\g{h}_{\g{p}}^\perp$, and using that
the elements of $\ad(\g{t})$ are skew-symmetric, together with
$\langle Z,Z\rangle=2$, yields $0=-b\lVert
X\rVert^4+2b\lVert X\rVert^4=b\lVert X\rVert^4$, which implies $X=0$.

Now take $W\in\g{w}$. Then $JW$ is orthogonal to
$\g{h}_{\g{p}}^\perp$ and thus there exists $R\in\g{t}$ such
that $R+JW\in\g{h}$. Since $\g{h}$ is orthogonal to
$[\g{h}_{\g{p}}^\perp,\g{h}_{\g{p}}^\perp]$, we obtain using~\eqref{E:bracket theta X, Z},
\[
0=\langle R+JW,[aB+b(1-\theta)Z,(1-\theta)W]\rangle
=\langle R+JW,\frac{a}{2}(1+\theta)W-b(1+\theta)JW\rangle
=-b\lVert W\rVert^2.
\]
Since $W\in\g{w}$ is arbitrary we conclude that $\g{w}=0$. Then
we can take $U\in\g{g}_\alpha$ and $S$, $T\in\g{t}$ such that
$S+U$, $T+JU\in\g{h}$. Since $\g{h}$ is orthogonal to
$\g{h}_{\g{p}}^\perp$, we get
\[
0=\langle [S+U,T+JU],aB+b(1-\theta)Z\rangle
=\langle [S,JU]-[T,U]+\frac{1}{2}\lVert U\rVert^2 Z,aB+b(1-\theta)Z\rangle
=b\lVert U\rVert^2,
\]
so $U=0$. This implies $\g{g}_\alpha=0$, which is not possible in $\C
H^n$, $n\geq 2$. Therefore $b=0$.
\end{proof}

Let us denote by $(\cdot)_{\g{a}\oplus\g{n}}$ and $(\cdot)_{\g{t}}$
the orthogonal projections onto $\g{a}\oplus\g{n}$ and $\g{t}$,
respectively. The next step of our proof is:

\begin{proposition}\label{P:Adg}
There exists an isometry $g\in G$ such that $(\Ad(g)\g{h})_{\g{a}\oplus\g{n}}=\g{s}_{V,\g{v}}$ for some subspace
$V$ of $\g{a}$ and some real subspace $\g{v}$ of $\g{g}_\alpha$.
\end{proposition}

\begin{proof}
Recall that $\g{h}_{\g{p}}^\perp=\R(aB+(1-\theta)X)\oplus(1-\theta)\g{w}$, where $a\in\R$, $X\in\g{g}_\alpha\ominus\g{w}$, and $\R X\oplus\g{w}$ is a real subspace of $\g{g}_\alpha$.

If $a=0$ the conclusion is obvious with $g=1$, the identity isometry in $G$, $V=\g{a}$ and $\g{v}=\R X+\g{w}$. So assume $a\neq 0$ and renormalize $X$ so that $a=1$. If $X=0$ the result is again obvious just by taking $g=1$ and $\g{v}=\g{w}$, so we will also assume $X\neq 0$.

Under these circumstances, we can define the isometry
\[
g=\Exp\left(-\frac{2}{\lVert X\rVert^2}X\right).
\]
We will show that
$(\Ad(g)\g{h})_{\g{a}\oplus\g{n}}=\g{s}_{\g{a},\g{w}}$. Note
that $\g{h}\cap\g{k}=\g{h}\cap\g{t}$ and
$\Ad(g)(\g{h})\cap\g{k}=\Ad(g)(\g{h})\cap\g{t}$ are conjugate
because $H$ and $I_g(H)$ are conjugate and all their orbits are
principal. Since
$\dim(\Ad(g)\g{h})_{\g{a}\oplus\g{n}}=\dim\Ad(g)\g{h}
-\dim(\Ad(g)(\g{h})\cap\g{t})=\dim\g{h}-\dim(\g{h}\cap\g{t})
=\dim\g{h}_{\g{a}\oplus\g{n}}=\dim\g{s}_{\g{a},\g{w}}$, it
suffices to show that
$(\Ad(g)\g{h})_{\g{a}\oplus\g{n}}\subset\g{s}_{\g{a},\g{w}}$.

From the definition of $\g{h}_{\g{p}}^\perp$ it follows that
\[\g{h}_{\g{a}\oplus\g{n}}=\R(\lVert X\rVert^2
B-X)\oplus(\g{g}_\alpha\ominus(\R
X\oplus\g{w}))\oplus\g{g}_{2\alpha}.
\]

Assume first that $T\in\g{h}\cap\g{t}$. If
$U\in\g{g}_\alpha\ominus(\R X\oplus\g{w})$ then there exists
$T_U\in\g{t}$ such that $T_U+U\in\g{h}$. Thus,
$[T,U]=[T,T_U+U]\in\g{h}$, and hence
$0=\langle[T,U],B+(1-\theta)X\rangle=-\langle[T,X],U\rangle$. This, together with the skew-symmetry of $\ad(T)$, implies $[T,X]\in\g{w}$. On the other hand, since
$\g{h}$ is orthogonal to $[\g{h}_{\g{p}}^\perp,\g{h}_{\g{p}}^\perp]$,
for any $W\in\g{w}$ we have, by~\eqref{E:auxiliary},
\[0=\langle
T,[B+(1-\theta)X,(1-\theta)W]\rangle=\langle T,\frac{1}{2}(1+\theta)W
-(1+\theta)[\theta
X,W]\rangle=-2\langle[T,X],W\rangle.\] Therefore $[T,X]=0$, and
$\Ad(g)$ acts as the identity on $\g{h}\cap\g{t}$, so
$(\Ad(g)(\g{h}\cap\g{t}))_{\g{a}\oplus\g{n}}=0$.

Let $S\in\g{t}$ be such that $S+\lVert X\rVert^2 B-X\in\g{h}$. For
any $U\in\g{g}_\alpha\ominus(\R X\oplus\g{w})$ there exists $T_U\in\g{t}$ such
that $T_U+U\in\g{h}$. Then, $0=\langle[S+\lVert X\rVert^2
B-X,T_U+U],B+(1-\theta)X\rangle=\langle[S,U],X\rangle$, so
$[S,X]\in\g{w}$. On the other hand, by~\eqref{E:auxiliary},
\[0=\langle S+\lVert X\rVert^2
B-X,[B+(1-\theta)X,(1-\theta)W]\rangle=-\langle S,(1+\theta)[\theta
X,W]\rangle=-2\langle[S,X],W\rangle\] for any $W\in\g{w}$. Thus,
$[S,X]=0$. This implies
\[
\Ad(g)(S+\lVert X\rVert^2 B-X)=S+\lVert
X\rVert^2 B-X-\frac{2}{\lVert X\rVert^2}\left(-\frac{\lVert
X\rVert^2}{2}X\right)=S+\lVert X\rVert^2 B,
\]
whose projection onto
$\g{a}\oplus\g{n}$ is $\lVert X\rVert^2 B\in\g{s}_{\g{a},\g{w}}$.

Now let $U\in\g{g}_\alpha\ominus(\R X\oplus\g{w})$ and $T_U\in\g{t}$
such that $T_U+U\in\g{h}$. Clearly, $[X,U]$,
$[X,[X,T_U]]\in\g{g}_{2\alpha}$. On the other hand, using~\eqref{E:auxiliary}, \[0=\langle
T_U+U,[B+(1-\theta)X,(1-\theta)W]\rangle=-\langle
T_U,(1+\theta)[\theta X,W]\rangle=-2\langle [T_U,X],W\rangle\] for any
$W\in\g{w}$. This and the skew-symmetry of $\ad(T_U)$ imply $[X,T_U]\in\g{g}_\alpha\ominus(\R
X\oplus\g{w})$, and thus
\begin{align*}
&(\Ad(g)(T_U+U))_{\g{a}\oplus\g{n}}
=\left(T_U+U-\frac{2}{\lVert
X\rVert^2}[X,T_U]-\frac{2}{\lVert
X\rVert^2}[X,U]+\frac{2}{\lVert
X\rVert^4}[X,[X,T_U]]\right)_{\g{a}\oplus\g{n}}\\
&\quad=U-\frac{2}{\lVert
X\rVert^2}[X,T_U]-\frac{2}{\lVert
X\rVert^2}[X,U]+\frac{2}{\lVert X\rVert^4}[X,[X,T_U]]
\in\g{g}_\alpha\ominus(\R X\oplus\g{w})\oplus\g{g}_{2\alpha}
\subset\g{s}_{\g{a},\g{w}}.
\end{align*}

Finally, let $T_Z\in\g{t}$ such that $T_Z+Z\in\g{h}$. Using~\eqref{E:auxiliary} we get \[0=\langle
T_Z+Z,[B+(1-\theta)X,(1-\theta)W]\rangle=-2\langle[T_Z,X],W\rangle\]
for any $W\in\g{w}$, thus showing that
$[T_Z,X]\in\g{g}_\alpha\ominus(\R X\oplus\g{w})$. Since
$[X,[X,T_Z]]\in\g{g}_{2\alpha}$ we get
\begin{align*}
&(\Ad(g)(T_Z+Z))_{\g{a}\oplus\g{n}}
=\left(T_Z+Z-\frac{2}{\lVert
X\rVert^2}[X,T_Z]+\frac{2}{\lVert
X\rVert^4}[X,[X,T_Z]]\right)_{\g{a}\oplus\g{n}}\\
&\quad=Z-\frac{2}{\lVert
X\rVert^2}[X,T_Z]+\frac{2}{\lVert X\rVert^4}[X,[X,T_Z]]
\in\g{g}_\alpha\ominus(\R X\oplus\g{w})\oplus\g{g}_{2\alpha}
\subset\g{s}_{\g{a},\g{w}}.
\end{align*}

Altogether this implies
$(\Ad(g)\g{h})_{\g{a}\oplus\g{n}}\subset\g{s}_{\g{a},\g{w}}$ and
since the dimensions are the same, equality follows as we had
claimed.
\end{proof}

So far we have shown that a group $H$ which induces a polar foliation
of $\C H^n$, $n\geq 2$, has a Lie subalgebra $\g{h}$ that can be
assumed to be contained in a maximally noncompact Borel subalgebra
$\g{t}\oplus\g{a}\oplus\g{n}$ in such a way that its projection onto
$\g{a}\oplus\g{n}$ satisfies
$\g{h}_{\g{a}\oplus\g{n}}=\g{s}_{V,\g{w}}$, where $\g{s}_{V,\g{w}}$
is one of our model examples with $V\subset\g{a}$ a linear subspace,
and $\g{w}\subset\g{g}_\alpha$ a real subspace. Our aim in what
follows is to show that the actions of $H$ and $S_{V,\g{w}}$ are
orbit equivalent.
The fundamental part of the proof of this fact is contained in the following assertion:

\begin{proposition}\label{P:normalizes}
With the notations as above, let $\g{h}_{\g{t}}$ denote the orthogonal projection of $\g{h}$ onto $\g{t}$. Then $[\g{h}_{\g{t}},\g{w}]=0$. In particular, $\g{h}_{\g{t}}$ normalizes $V\oplus(\g{n}\ominus\g{w})$.
\end{proposition}

\begin{proof}
We begin with some general observations regarding $\g{h}$. Let $W_1$, $W_2\in\g{w}$, $U\in \g{s}_{V,\g{w}}=V\oplus(\g{g}_\alpha\ominus\g{w})\oplus\g{g}_{2\alpha}$, and $T\in\g{t}$ such that $T+U\in\g{h}$. Since $\g{h}$ is orthogonal to $[\g{h}_{\g{p}}^\perp,\g{h}_{\g{p}}^\perp]$ we have by~\eqref{E:auxiliary},
\[0=\langle T+U,[(1-\theta)W_1,(1-\theta)W_2]\rangle=-\langle T,(1+\theta)[\theta W_1,W_2]\rangle=-2\langle [T,W_1],W_2\rangle.\]
This already implies $[\g{h}_{\g{t}},\g{w}]\subset\g{g}_\alpha\ominus\g{w}$.

Now let $S\in\g{h}\cap\g{t}$. Take $X\in\g{g}_\alpha\ominus\g{w}$ and $T\in\g{t}$ such that $T+X\in\g{h}$. Then $[S,X]=[S,T+X]\in\g{h}$ and thus, for any $W\in\g{w}$ we get $0=\langle [S,X],(1-\theta)W\rangle=-\langle [S,W],X\rangle$. Together with $[\g{h}_{\g{t}},\g{w}]\subset\g{g}_\alpha\ominus\g{w}$, this yields $[\g{h}\cap\g{t},\g{w}]=0$.

For $W\in\g{w}$ we define the map $F_W\colon\g{s}_{V,\g{w}}\to\g{g}_\alpha\ominus\g{w}$, $U\mapsto[W,T_U]$, where $T_U\in\g{t}$ is any vector satisfying $T_U+U\in\g{h}$. This map is well-defined because given $S_U$, $T_U\in\g{t}$ such that $S_U+U$, $T_U+U\in\g{h}$, we have $S_U-T_U=(S_U+U)-(T_U+U)\in\g{h}\cap\g{t}$, and from $[\g{h}\cap\g{t},\g{w}]=0$ it follows that $[W,S_U]=[W,T_U]$. The map is clearly linear, and its image is contained in $\g{g}_\alpha\ominus\g{w}$ because $[\g{h}_{\g{t}},\g{w}]\subset\g{g}_\alpha\ominus\g{w}$. The assertion will follow if we show that $F_W=0$ for any $W\in\g{w}$.

Let $H\in V$, $X\in\g{g}_\alpha\ominus\g{w}$, $b\in\R$, and $S$, $T\in\g{t}$ such that $T+H+bZ$, $S+X\in\g{h}$. Since $[\g{t},\g{a}\oplus\g{g}_{2\alpha}]=0$, we get $[T,X]+\frac{1}{2}\langle H,B\rangle X=[T+H+bZ,S+X]\in\g{h}$. Thus, for any $W\in\g{w}$ we have $0=\langle [T,X]+\frac{1}{2}\langle H,B\rangle X,(1-\theta)W\rangle=\langle [W,T],X\rangle$. Since $X\in\g{g}_\alpha\ominus\g{w}$ is arbitrary and $[\g{h}_{\g{t}},\g{w}]\subset\g{g}_\alpha\ominus\g{w}$, we get $F_W(V\oplus\g{g}_{2\alpha})=0$ for any $W\in\g{w}$.

Let us denote by $\tilde{F}_W\colon \g{g}_\alpha\ominus\g{w}\to\g{g}_\alpha\ominus\g{w}$ the restriction of $F_W$ to $\g{g}_\alpha\ominus\g{w}$. It now suffices to prove that $\tilde{F}_W=0$ for any $W\in\g{w}$.

Let $X$, $Y\in\g{g}_\alpha\ominus\g{w}$ and $T_X$, $T_Y\in\g{t}$ be such that $T_X+X$, $T_Y+Y\in\g{h}$. Then $[T_X,Y]-[T_Y,X]+\frac{1}{2}\langle JX,Y\rangle Z=[T_X+X,T_Y+Y]\in\g{h}$, and thus, for any $W\in\g{w}$, taking inner product with $(1-\theta)W$, we get $\langle [W,T_X],Y\rangle=\langle [W,T_Y],X\rangle$. This readily implies $\langle\tilde{F}_W(X),Y\rangle=\langle\tilde{F}_W(Y),X\rangle$ for any $X$, $Y\in\g{g}_\alpha\ominus\g{w}$, so $\tilde{F}_W$ is a self-adjoint endomorphism of $\g{g}_\alpha\ominus\g{w}$, and hence, it is diagonalizable with real eigenvalues.

Assume, by contradiction, that there exists $W\in\g{w}$ such that $\tilde{F}_W\neq 0$. Then there is and eigenvalue $\lambda\neq 0$ of $\tilde{F}_W$ and a nonzero vector $X\in\g{g}_\alpha\ominus\g{w}$ such that $\tilde{F}_W(X)=\lambda X$. Let $T_X\in\g{t}$ be such that $T_X+X\in\g{h}$; then $[W,T_X]\neq 0$ because $\tilde{F}_W(X)\neq 0$.

We define $g=\Exp(-\frac{1}{\lambda}W)$. The groups $H$ and
$I_g(H)$ are conjugate, and since $H$ induces a foliation so
does $I_g(H)$. Since all the orbits of a homogeneous foliation
are principal and $H$ and $I_g(H)$ are conjugate, so are their
isotropy groups at the origin, $H_o$ and $I_g(H)_o$ respectively. In
particular, their Lie algebras must have the same dimension, that
is, $\dim(\g{h}\cap\g{k})=\dim(\Ad(g)(\g{h})\cap\g{k})$. We will see that this is not true, hence obtaining our contradiction.

Clearly, $\g{h}\cap\g{k}=\g{h}\cap\g{t}$. On the other hand, since $[\g{h}\cap\g{t},\g{w}]=0$, we get $\Ad(g)\vert(\g{h}\cap\g{t})=1_{\g{h}\cap\g{t}}$, the identity on $\g{h}\cap\g{t}$, so $\g{h}\cap\g{t}\subset\Ad(g)(\g{h})\cap\g{k}$. Finally, let $b=\frac{1}{2\lambda}\langle JW,X\rangle-\frac{1}{4\lambda^2}\langle JW,[W,T_X]\rangle$ and $S\in\g{t}$ such that $S+bZ\in\g{h}$. Recall that $[W,S]=0$ because of $F_W(\g{g}_{2\alpha})=0$. Then,
\begin{align*}
&\Ad(g)(T_X+X+S+bZ)
=\sum_{k=0}^\infty\frac{1}{k!}
\ad\!\left(-\frac{1}{\lambda}W\right)^k(T_X+X+S+bZ)\\
&\qquad=T_X+X+S+bZ-\frac{1}{\lambda}([W,T_X]+[W,X])
+\frac{1}{2\lambda^2}[W,[W,T_X]]\\
&\qquad=T_X+S+X-\frac{1}{\lambda}\tilde{F}_W(X)
+\left(b-\frac{1}{2\lambda}\langle JW,X\rangle+\frac{1}{4\lambda^2}\langle JW,[W,T_X]\rangle\right)Z\\
&\qquad=T_X+S\in\g{k}.
\end{align*}
Now, $T_X+S\notin\g{t}\cap\g{h}$ as $[W,T_X]\neq 0$, $[S,\g{w}]=0$, and $[\g{h}\cap\g{t},\g{w}]=0$. This exhibits the fact that $\g{h}\cap\g{t}$ is strictly contained in $\Ad(g)(\g{h})\cap\g{k}$, which contradicts the property that $H$ induces a foliation. Therefore $\tilde{F}_W=0$ for all $W\in\g{w}$, and we get $[\g{h}_{\g{t}},\g{w}]=0$.

In order to obtain the last assertion, first observe that $\g{t}$ centralizes $\g{a}\oplus\g{g}_{2\alpha}$. Now if $X\in\g{g}_\alpha$, $T\in\g{h}_{\g{t}}$, and $W\in\g{w}$, we get $\langle[T,X],W\rangle=-\langle[T,W],X\rangle=0$, from where it follows that $\g{h}_{\g{t}}$ normalizes $\g{g}_\alpha\ominus\g{w}$.
\end{proof}

Now we are ready for the final step of the proof.

\begin{proposition}
Let $\g{t}\oplus\g{a}\oplus\g{n}$ be a maximally noncompact Borel
subalgebra of $\g{g}$. Let $\g{h}$ be a subalgebra of
$\g{t}\oplus\g{a}\oplus\g{n}$ such that
$\g{h}_{\g{a}\oplus\g{n}}=\g{s}_{V,\g{w}}$. Assume that the
orbits of the connected subgroup $H$ of $G$ whose Lie algebra is
$\g{h}$ form a homogeneous polar foliation of $\C H^n$. Then the actions of $H$
and $S_{V,\g{w}}$ on $\C H^n$ are orbit equivalent.
\end{proposition}

\begin{proof}
First recall from Proposition~\ref{P:normalizes} that $\g{h}_{\g{t}}$ normalizes
$V\oplus(\g{g}_\alpha\ominus\g{w})\oplus\g{g}_{2\alpha}
=V\oplus(\g{n}\ominus\g{w})$. We also denote by $H_{\g{t}}$ the connected subgroup of $SU(1,n)$ whose Lie algebra is $\g{h}_{\g{t}}$.

We now show that an element $g\in H$ can be written as $g=tan$ where
$t\in H_{\g{t}}$, $a\in\Exp(V)\subset A$ and $n\in S_{\{0\},\g{w}}$,
where as usual $S_{\{0\},\g{w}}$ is the connected subgroup of $AN$  whose
Lie algebra is $\g{n}\ominus\g{w}$. First write $g=\Exp(T+aB+X)$ with
$T\in\g{t}$, $a\in\R$, and $X\in\g{n}\ominus\g{w}$. By the Iwasawa
decomposition we can write $g=\Exp(S)\Exp(bB)\Exp(Y)$, with $S\in\g{t}$,
$b\in\R$, and $Y\in\g{n}$. The Baker-Campbell-Hausdorff formula
yields
\begin{align*}
\Exp(Y) &=\Exp(-bB)\Exp(-S)\Exp(T+aB+X)=\Exp(-S-bB)\Exp(T+aB+X)\\
        &=\Exp(T-S+(a-b)B+X-[S,X]-b[B,X]+\dots),
\end{align*}
where the dots involve linear combinations of iterated brackets of
$X$ with $S$ and $B$. Comparing both sides of the equation we
immediately get $S=T$ and $b=a$. Now, since $\g{h}_\g{t}\oplus\g{a}$
normalizes $\g{n}\ominus\g{w}$, the
right-hand side of the equation above is in $\g{n}\ominus\g{w}$, and
so is $Y$, as we wanted to show.

From Proposition~\ref{P:normalizes} we also have $[\g{h}_{\g{t}},\g{w}]=0$,
which obviously implies $\Ad(\Exp(T))\xi=\xi$ for any
$\xi\in(1-\theta)\g{w}$ and $T\in\g{h}_{\g{t}}$. Let $g\in H$ be
arbitrary and write $g^{-1}=tan$ with $t\in H_{\g{t}}\subset K$,
$a\in \Exp(V)\subset A$ and $n\in S_{\{0\},\g{w}}$ as explained above.
Then $g=n^{-1}a^{-1}t^{-1}$, and since $A$ normalizes $S_{\{0\},\g{w}}$
we can write $n^{-1}a^{-1}=a^{-1}n'$ with $n'\in S_{\{0\},\g{w}}$. Thus,
$g(o)=a^{-1}n't^{-1}(o)=a^{-1}n'(o)$ and hence $H\cdot o\subset
S_{V,\g{w}}\cdot o$. Since both orbits $H\cdot o$ and
$S_{V,\g{w}}\cdot o$ have the same dimension and are connected and
complete we conclude $H\cdot o=S_{V,\g{w}}\cdot o$.

Now, let $p=\exp_o(\xi)$ with $\xi\in\nu_o(H\cdot
o)\cong(1-\theta)\g{w}\subset\g{p}$. Using the fact that $H$ acts
isometrically on $\C H^n$, $n\geq 2$, and that
$t^{-1}_*\xi=\Ad(t^{-1})\xi=\xi$ we get
\[
g(p)=\exp_{g(o)}(g_*\xi)=\exp_{an't^{-1}(o)}((an't^{-1})_*\xi)
=\exp_{an'(o)}((an')_*\xi)=an'(p).
\]
Hence, $H\cdot p\subset S_{V,\g{w}}\cdot p$, and thus equality
follows. Since the action of $H$ is polar, all the orbits can be
obtained in this way, and so, $H$ and $S_{V,\g{w}}$ have the same
orbits.
\end{proof}

This concludes the proof of the Main Theorem.


\end{document}